\documentclass{svproc}
\usepackage{amssymb}
\usepackage{mathrsfs}
\usepackage{graphicx}
\usepackage{subfigure}

\usepackage{url}

\newtheorem{observation}[theorem]{Observation}
\newcommand{\ERCagreement}{{\begin{minipage}{.64\textwidth}This paper is part of a project that has received funding from the European Research Council (ERC) under the European Union's Horizon 2020 research and innovation programme (grant agreement No 810115 -- {\sc Dynasnet}). \end{minipage}\hfill\begin{minipage}{.32\textwidth}\includegraphics[width=\textwidth]{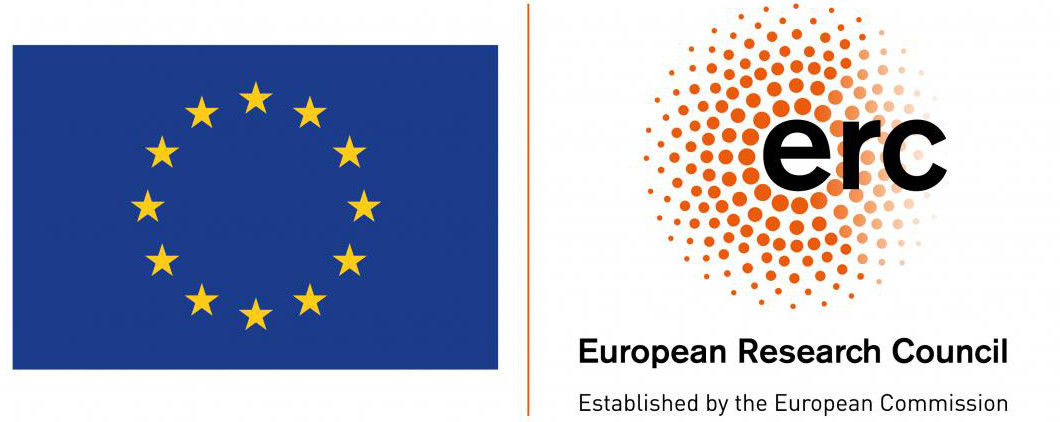}\end{minipage}\hfill}}

\begin{document}
\mainmatter              
\title{On asymmetric hypergraphs\thanks\ERCagreement}
\titlerunning{On asymmetric hypergraphs}  
%
\author{Yiting Jiang\inst{1} \and Jaroslav Ne\v{s}et\v{r}il\inst{2}}
\authorrunning{Yiting Jiang, and Jaroslav Ne\v{s}et\v{r}il} 
%
\tocauthor{Yiting Jiang, and Jaroslav Ne\v{s}et\v{r}il}
\institute{Universit\'e de Paris, CNRS, IRIF, F-75006, Paris, France\\and Department of Mathematics, Zhejiang Normal University, China
    \email{yjiang@irif.fr},\\ 
	\and
	Computer Science Institute of Charles University (IUUK and ITI)\\
	Malostransk\' e n\' am.25, 11800 Praha 1, Czech Republic\\
	\email{nesetril@iuuk.mff.cuni.cz}}
	
\maketitle              

\begin{abstract}
In this paper, we prove that for any $k\ge 3$, there exist infinitely many minimal asymmetric $k$-uniform hypergraphs. This is in a striking contrast to $k=2$, where it has been proved recently that there are exactly $18$ minimal asymmetric graphs.

We also determine, for every $k\ge 1$, the minimum size of an asymmetric $k$-uniform hypergraph.
\end{abstract}

\section{Introduction}

Let us start with graphs: A graph $G$ is called \emph{asymmetric} if it does not have a non-identical automorphism. Any non-asymmetric graph is also called \emph{symmetric} graph. A graph $G$ is called \emph{minimal asymmetric} if $G$ is asymmetric and every non-trivial induced subgraph of $G$ is symmetric. (Here $G'$ is a \emph{non-trivial subgraph} of $G$ if $G'$ is a subgraph of $G$ and $1<|V(G')|<|V(G)|$.) In this paper all graphs are finite.

It is a folkloristic result that most graphs are asymmetric. In fact, as shown by Erd\H{o}s and R\'{e}nyi \cite{AGER} most graphs on large sets are asymmetric in a very strong sense. 
The paper \cite{AGER} contains many extremal results (and problems), which motivated further research on extremal properties of asymmetric graphs, see e.g. \cite{NS}, \cite{SS}.
This has been also studied in the context of the reconstruction conjecture \cite{CTAT}, \cite{VM}.

The second author bravely conjectured long time ago that there are only finitely many minimal asymmetric graphs, see e.g. \cite{BOOK}. Partial results were given in \cite{NS92}, \cite{Sab91}, \cite{Woj96}. Recently this conjecture has been confirmed by Pascal Schweitzer and Patrick Schweitzer \cite{MSH} (the list of $18$ minimal asymmetric graphs has been isolated in \cite{NS92}):

\begin{theorem}\cite{MSH}
	\label{thm1}
	There are exactly 18 minimal asymmetric undirected graphs up to isomorphism.
\end{theorem}

An {\emph{involution}} of a graph $G$ is any non-identical automorphism $\phi$ for which $\phi\circ\phi$ is an identity. It has been proved in \cite{MSH} that all minimal asymmetric graphs are in fact minimal involution-free graphs.

\medskip
In this paper, we consider analogous questions for \emph{$k$-graphs} (or $k$-uniform hypergraphs), i.e. pairs $(X,\mathscr{M})$ where $\mathscr{M}\subseteq {X\choose k}=\{A\subseteq X;|A|=k\}$. 
Induced subhypergraphs, asymmetric hypergraphs and minimal asymmetric hypergraphs are defined analogously as for graphs. Instead of hypergraphs we often speak just about $k$-graphs.

We prove two results related to minimal asymmetric $k$-graphs.

Denote by $n(k)$ the minimum number of vertices of an asymmetric $k$-graph.

\begin{theorem}
	\label{main1}
	$n(2)=6$, $n(3)=k+3$, $n(k)=k+2$ for $k\ge 4$.	
\end{theorem}

Theorem \ref{main1} implies the existence of small minimal asymmetric $k$-graphs. Our second result disproves analogous minimality conjecture (i.e. a result analogous to Theorem \ref{thm1}) for $k$-graphs.

\begin{theorem}
	\label{main2}
	For every integer $k\ge 3$, there exist infinitely many $k$-graphs  that are minimal asymmetric.
\end{theorem}

In fact we prove the following stronger statement.

\begin{theorem}
	\label{main3}
	For every integer $k\ge 3$, there exist infinitely many $k$-graphs $(X, \mathscr{M})$ such that 
	\begin{itemize}
		\item[1.] $(X, \mathscr{M})$ is asymmetric. 
		\item[2.] If $(X', \mathscr{M}')$ is a $k$-subgraph of $(X, \mathscr{M})$ with at least two vertices, then $(X', \mathscr{M}')$ is symmetric.
	\end{itemize}
\end{theorem}

Such $k$-graphs we call \emph{strongly minimal asymmetric}. So strongly minimal asymmetric $k$-graphs do not contain any non-trivial
(not necessarily induced) asymmetric $k$-subgraph. Note that some of the minimal asymmetric graphs fail to be strongly minimal.

\medskip

Theorem \ref{main3} is proved by constructing a sequence of strongly minimal asymmetric $k$-graphs. We have two different constructions of increasing strength: In section 3 we give a construction with all vertex degrees bounded by $3$. A stronger construction which yields minimal asymmetric $k$-graphs ($k\ge 6$) with respect to involutions is presented due to space limitations in the Appendix A.

\section{The proof of Theorem \ref{main1}}

\begin{lemma}
	\label{lem}
	For $k\ge 3$, we have $n(k)\ge k+2$. 
\end{lemma}

\begin{proof}
	Assume that there exists an asymmetric $k$-graph $(X, \mathscr{M})$ with $|X|=k+1$.
	If for each vertex $u\in X$, there is a hyperedge $M\in\mathscr{M}$ such that $u\notin M$, then $\mathscr{M}={X\choose k}$, which  is symmetric.
	Otherwise there exists $u,v\in X$ such that $\{u,v\}\subset M$  for every edge $M\in\mathscr{M}$, or there exist $u',v'\in X$ and $M_1,M_2\in\mathscr{M}$ such that $u'\notin M_1$ and $v'\notin M_2$.
	In the former case, there is an automorphism $\phi$ of $(X, \mathscr{M})$ such that $\phi(u)=v$ and $\phi(v)=u$. In the later case there is an automorphism $\phi$ of $(X,\mathscr{M})$ such that $\phi(u')=v'$ and $\phi(v')=u'$. In both cases we have a contradiction. \qed
\end{proof}

For a $k$-graph $G=(X,\mathscr{M})$, the \emph{set-complement} of $G$ is defined as a $(|X|-k)$-graph $\bar{G}=(X,\bar{\mathscr{M}})=(X,\{X-M|M\in\mathscr{M}\})$.
Denote by $Aut(G)$ the set of all the automorphisms of $G$ and thus we have $Aut(G)=Aut(\bar{G})$. 
We define the degree of a vertex $v$ in a $k$-graph $G$ as 
$d_G(v)=|\{M\in\mathscr{M};v\in M\}|$. 

\begin{lemma}
	\label{lem1}
	For $k\ge 4$, we have $n(k)=k+2$.
\end{lemma}


\begin{figure}[htbp]
		\centering
		\subfigure[$X_1$]{
		\includegraphics[scale=0.5]{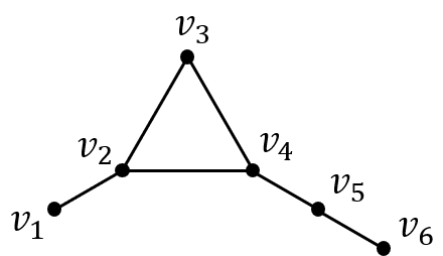}
        }
    \qquad\qquad\qquad\qquad
		\subfigure[$T_{k+2}$]{
		\centering
		\includegraphics[scale=0.5]{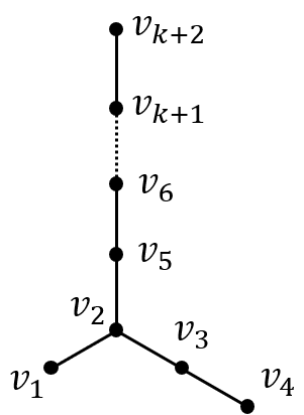}
       }
   \caption{}
\end{figure}


\begin{proof}
	First, we construct an asymmetric $k$-graph $(X,\mathscr{M})$ with $|X|=k+2$ for each $k\ge 4$.
	Examples of such graphs $X_1$ and $T_{k+2}$ are depicted in Fig. 1.

	For $k=4$, take the set-complement of $X_1$.
	For every $k\ge 5$, take the set-complement of $T_{k+2}$.  
	It is known that $X_1$ and $T_{k+2}$ ($k\ge 5$) are asymmetric.
	Thus set-complements $\bar{X}_1$ and $\bar{T}_{k+2}$ ($k\ge 5$) are also asymmetric $k$-graphs. 
	
	Let $H$ be any of these set-complements.
	$H$ is minimum as any non-trivial induced sub-$k$-graph $H'\subset H$ has no more than $k+1$ elements. Thus we can use Lemma \ref{lem}. \qed
\end{proof}

The case $k=3$ (i.e. proof of $n(3)=6$) is a (rather lenghthy) case analysis which has to be ommitted here.

\section{Proofs of Theorem \ref{main3}}

In this section, we outline two different proofs of Theorem \ref{main3}.

Firstly, we give a proof with bounded degrees.

For $k\ge 3$, $t\ge k-2$, we define the following $k$-graphs.

\medskip
$G_{k,t}=(X_{k,t},\mathscr{E}_{k,t})$, 

$X_{k,t}=\{v_i;i\in [tk]\}\cup\{u_i;i\in [tk]\}\cup\{v^{j}_i;i\in [tk], j\in [k-3]\}\}$,

$\mathscr{E}_{k,t}=\{E_i;i\in [tk]\}\cup\{E_{i,j};j\in[k-3],i=j+sk,s\in\{0,1,2,\cdots,t-1\}\}$, where $E_i=\{v_i,u_i,v^1_i,v^2_i,\cdots,v^{k-3}_i,v_{i+1}\}$, $E_{i,j}=\{v^j_i,v^j_{i+1},\cdots,v^j_{i+k-1}\}$ and  $y_{tk+t}=y_t$ for every $y\in\{v,u,v^1,v^2,\cdots,v^{k-3}\}$.

\medskip
$G^{\circ}_{k,t}=\{X_{k,t}\cup\{x\}, \mathscr{E}_{k,t}\cup\{E\}\}$, where $E=\{v_1,u_1,v^1_1,v^2_1,\cdots,v^{k-3}_1,x\}$.

\medskip
The $k$-graphs $G_{k,t}$ and $G^{\circ}_{k,t}$ are schematically depicted in Fig. 3.


\begin{figure}[htbp]
	\centering
	\subfigure[$G_{k,t}$]{
		\includegraphics[scale=0.2]{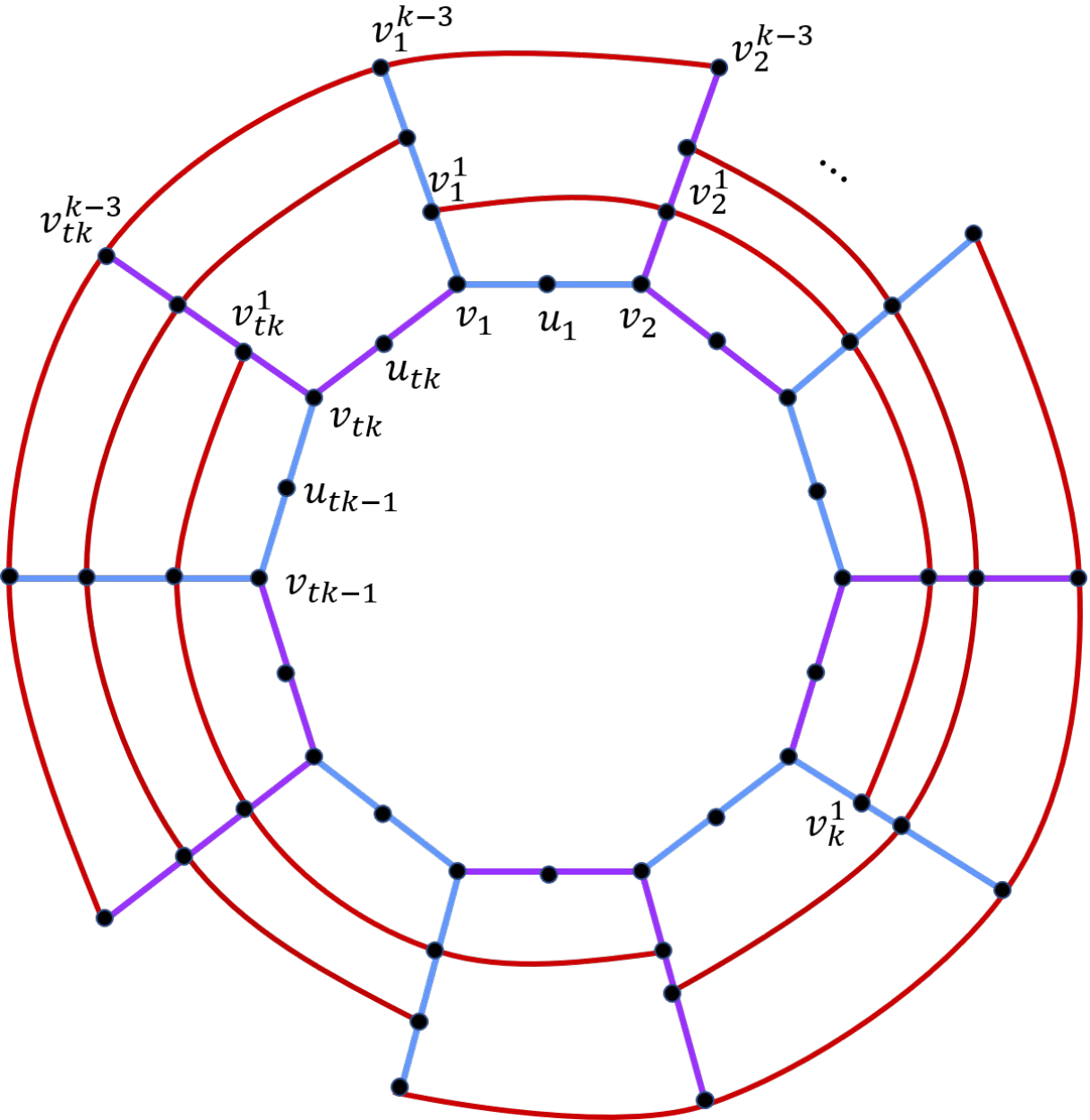}
	}
\qquad\qquad\qquad
	\subfigure[$G^{\circ}_{k,t}$]{
		\centering
		\includegraphics[scale=0.2]{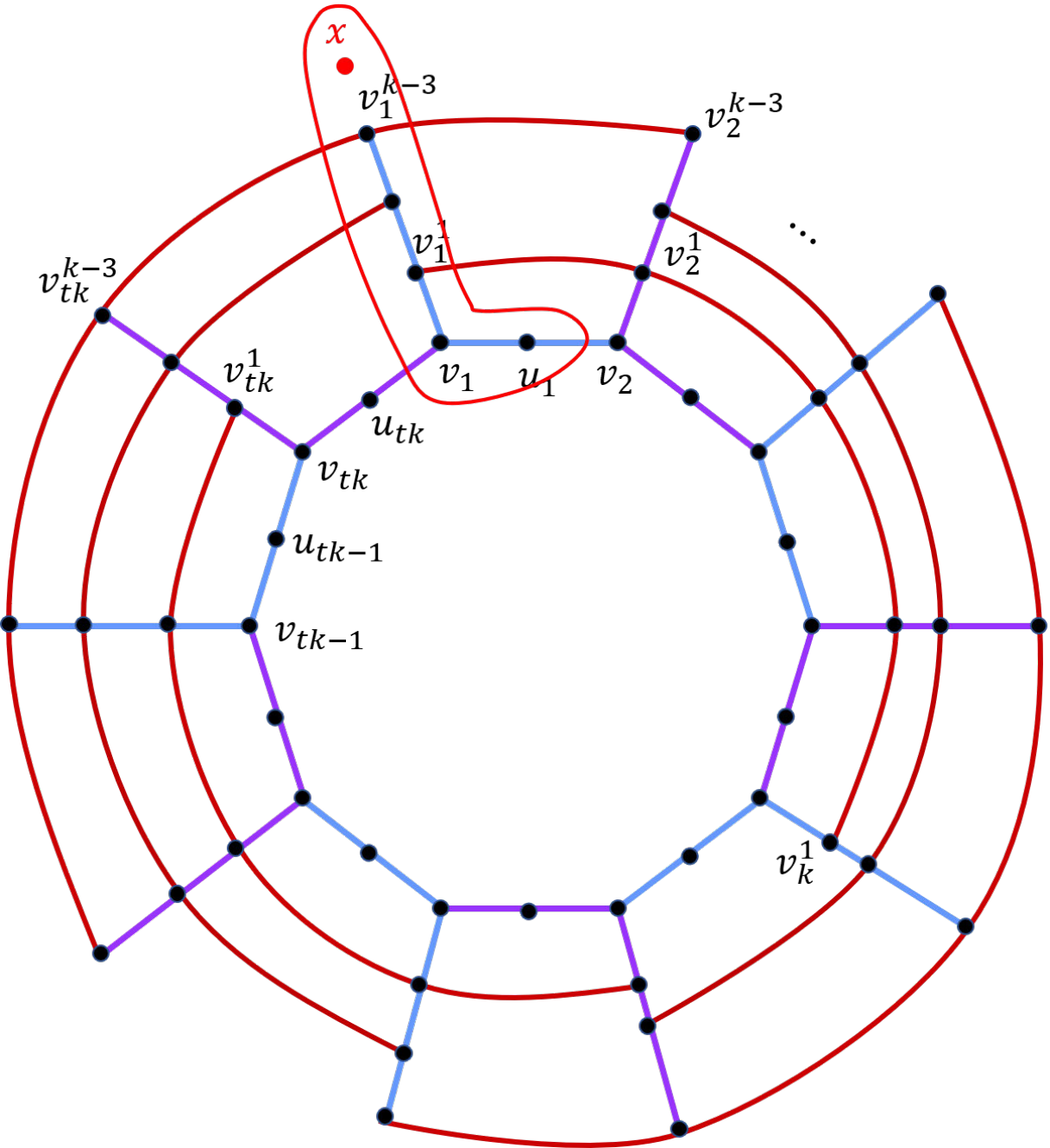}
	}
	\caption{}
\end{figure}


\medskip
The proof of Theorem \ref{main3} follows from the following two lemmas. The proofs are ommitted due to space limitations:

\begin{lemma}
	\label{lem7}
	\begin{itemize}
		\item[1)] The $k$-graph $G_{k,t}$ is symmetric and every non-identical automorphism $\phi$ of $G_k$ satisfies that there exists an integer $c$ (which does not divide $tk$) such that for every $i\in [tk]$, $j=i+c-\lfloor\frac{i+c}{tk}\rfloor$, $\phi(E_i)=E_j$ (i.e. for each vertex $v\in E_i$, $\phi(v)\in E_j$).
		\item[2)] The only automorphism of $G_{k,t}$ which leaves the set $E_1\setminus\{v_2\}$ invariant (i.e. for each vertex $v\in E_1\setminus\{v_2\}$, $\phi(v)\in E_1\setminus\{v_2\}$) is the identity. 
		\item[3)] Every non-trivial $k$-subgraph of $G_{k,t}$ containing the vertices in $E_1$ has a non-identical automorphism $\phi$ which leaves the set $E_1$ invariant.
	\end{itemize}
\end{lemma}

\begin{lemma}
	\begin{itemize}
		\item[1)] The $k$-graph $G^{\circ}_{k,t}$ is asymmetric. 
		\item[2)] Every non-trivial $k$-subgraph of $G^{\circ}_{k,t}$ has a non-identical automorphism.
	\end{itemize}
\end{lemma}

\medskip
It is easy to observe that the $k$-graphs $G^{\circ}_{k,t}$ have vertex degree at most $3$.
However note that in this construction, some of the strongly minimal asymmetric $k$-graphs $G^{\circ}_{k,t}$ are not minimal involution-free. In fact, when $k\ge 3$, $t\ge k-2$ is odd, the $k$-subgraph $G^{\circ}_{k,t}-x$ of $G^{\circ}_{k,t}$ is involution-free. 
However the most interesting form of Theorem \ref{main3} relates to minimal asymmetric graphs for involutions.
This can be stated as follows:
\begin{theorem}
	\label{main4}
	For every $k\ge 6$, there exist infinitely many $k$-graphs $(X, \mathscr{M})$ such that 
		\begin{itemize}
			\item[1.] $(X, \mathscr{M})$ is asymmetric. 
			\item[2.] If $(X', \mathscr{M}')$ is a $k$-subgraph of $(X, \mathscr{M})$ with at least two vertices, then $(X', \mathscr{M}')$ has an involution.
		\end{itemize}
\end{theorem}
This is more involved and the proof is attached as Appendix A.


\section{Concluding remarks}

Of course one can define the notion of asymmetric graph also for directed graphs, binary relations and $k$-nary relations $R\subseteq X^k$.

One has then the following analogy of Theorem 1:
There are exactly $19$ minimal asymmetric binary relations.
(These are symmetric orientations of $18$ minimal asymmetric (undirected) graphs and the single arc graph ($\{0,1\}$, $\{(0,1)\}$).)

Here is a companion problem about extremal asymmetric oriented graphs and one of the original motivation of the problem \cite{BOOK}:

Let $G=(V,E)$ be an asymmetric graph with at least two vertices. We say that $G$ is \emph{critical asymmetric} if for every $x\in V$ the graph $G-x=(V\setminus\{x\},\{e\in E;x\notin e\})$ fails to be asymmetric.
An oriented graph is a relation not containing two opposite arcs.

\smallskip
\noindent
{\bfseries{Conjecture 1}}
Let $G$ be an oriented asymmetric graph. Then it fails to be critical asymmetric. Explicitly: For every oriented asymmetric graph $G$, there exists $x\in V(G)$ such that $G-x$ is asymmetric. 

\smallskip
W\'{o}jcik \cite{Woj96} proved that a critical oriented asymmetric graph has to contain a directed cycle. 

This research indicates a particular role of binary structures with respect to automorphism and asymmetry. While for higher arities there are infinitely many minimal asymmetric graphs, for binary structures this may be always finite. We formulate this in graph language as follows:

Let $L$ be a finite set of colours. An $L$-graph is a finite graph where each edge gets one of the colours from $L$. Automorphisms are defined as colour preserving automorphisms. The following is a problem which generalizes the problem (and its solution [8]) which motivated the present note:

\smallskip
\noindent
{\bfseries{Conjecture 2}}
For every finite $L$, there are only finitely many minimal asymmetric $L$-graphs.

\smallskip
This is open  for any $|L| > 1$. By results of this paper we know that for non-binary languages the analogous problem has negative solution.
%
%

\clearpage

\section*{Appendix A}
\appendix
\setcounter{theorem}{0}
\setcounter{lemma}{0}
\setcounter{figure}{0}
\renewcommand{\thetheorem}{A\arabic{theorem}}
\renewcommand{\thelemma}{A\arabic{lemma}}
\renewcommand{\thefigure}{A\arabic{figure}}

We prove here Theorem \ref{main4} (recalled as follows):

\medskip
\noindent
{\bfseries{Theorem \ref{main4} }}\emph{For every $k\ge 6$, there exist infinitely many $k$-graphs $(X, \mathscr{M})$ such that 
	\begin{itemize}
		\item[1.] $(X, \mathscr{M})$ is asymmetric. 
		\item[2.] If $(X', \mathscr{M}')$ is a $k$-subgraph of $(X, \mathscr{M})$ with at least two vertices, then $(X', \mathscr{M}')$ has an involution.
\end{itemize}}

\medskip

To prove this theorem, we first construct the following $k$-graphs:

\medskip
$G_k=(X_k,\mathscr{M}_k)$, $X_k=\{v_1,v_2,\cdots,v_{2k-1}\}$, $\mathscr{M}_k=\{M_i=\{v_i,v_{i+1},\cdots,v_{i+k-1}\};i\in [k]\}$. 

\medskip
$G^*_k=(X^*_k,\mathscr{M}^*_k)$, $X^*_k=X_k\cup\{x\}$, $\mathscr{M}^*_k=\mathscr{M}_k\cup\{x,v_1,\cdots,v_{k-2},v_{k+2}\}$.

\medskip

These $k$-graphs are depicted in Figure A1 and A2.


\begin{figure}[!htp]
	\begin{center}
		\includegraphics[scale=0.3]{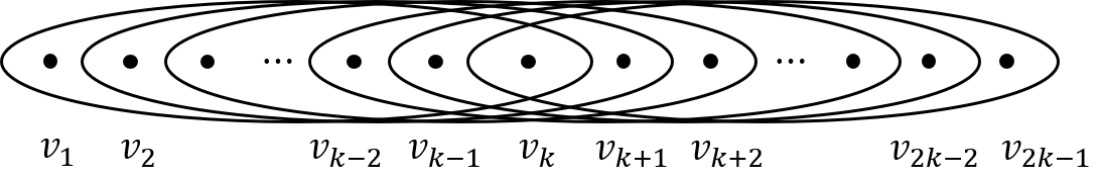}
	\end{center}
	\caption{The graph $G_k$}
\end{figure}


\begin{figure}[!htp]
	\begin{center}
		\includegraphics[scale=0.3]{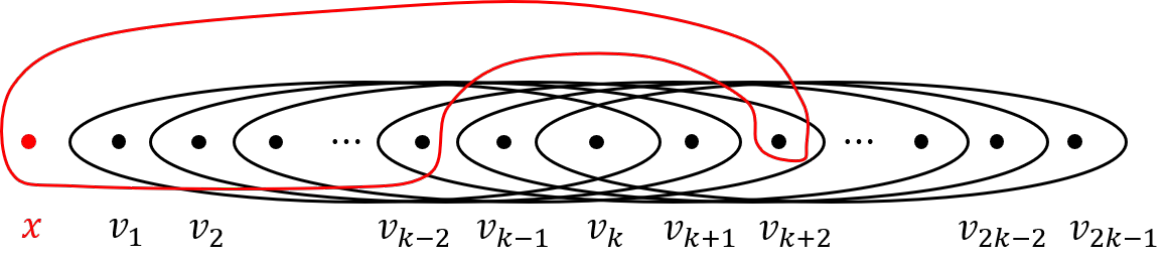}
	\end{center}
	\caption{The graph $G^*_k$}
\end{figure}


\begin{lemma}
	\label{lem3}
	\begin{itemize}
		\item[1)] The $k$-graph $G_k$ is symmetric and the only non-identical automorphism $\phi$ of $G_k$ satisfies that $\phi(v_i)=v_{2k-i}$ for every $i\in [2k-1]$.
		\item[2)] The only automorphism of $G_k$ which leaves the set $\{v_{2k-2}, v_{2k-1}\}$ invariant (i.e. $\{\phi(v_{2k-2}), \phi(v_{2k-1})\}=\{v_{2k-2}, v_{2k-1}\}$) is the identity. 
		\item[3)] Every non-trivial $k$-subgraph of $G_k$ containing vertices $v_{2k-2}$, $v_{2k-1}$ has an involution $\phi$ which leaves the set $\{v_{2k-2}, v_{2k-1}\}$ invariant.
	\end{itemize}
\end{lemma}

\begin{proof}
	The first property holds by observing the degree of each vertex in $G_k$. Then the second follows.
	
	To prove the third one, we assume that $G$ is a non-trivial $k$-subgraph of $G_k$ such that $G$ contains vertices $v_{2k-2}$, $v_{2k-1}$ and $j$ is the maximal index that $G$ contains the edge $M_j=\{v_j,v_{j+1},\cdots,v_{j+k-1}\}$. Let $i$ be the minimal index such that $G$ contains the edges $M_i, M_{i+1}, \dots, M_j$. Since $G$ is a nontrival $k$-subgraph of $G_k$, we have $j<k$ and $M_{j+1}$ is not in $G$ or $i>1$ and $M_{i-1}$ is not in $G$. It implies that $v_{i+k-2}$, $v_{i+k-1}$ share the same edges $M_i$, $M_{i+1}$, $\dots$, $M_j$. If $i\notin\{k-1,k\}$ then there is an involution $\phi$ of $G$ which leaves the set $\{v_{2k-2}, v_{2k-1}\}$ invariant, $\phi(v_{i+k-2})=v_{i+k-1}$ and $\phi(v_{i+k-1})=v_{i+k-2}$. If $i\in\{k-1,k\}$ and $G$ contains an edge $M_l$ ($1\le l<i-1$), then there is an involution $\phi$ of $G$ which leaves the set $\{v_{2k-2}, v_{2k-1}\}$ invariant, $\phi(v_{i-2})=v_{i-1}$ and $\phi(v_{i-1})=v_{i-2}$, as $M_l$ contains the vertex $v_{i-1}$ and $M_{i-1}$ is not in $G$. Now the remaining case is the edge set of $G$ is contained in $\{M_{k-1},M_k\}$, which is easy to observe that there is an involution $\phi$ of $G$ which leaves the set $\{v_{2k-2}, v_{2k-1}\}$ invariant. \qed
	
\end{proof}

\begin{lemma}
	\label{lem4}
	\begin{itemize}
		\item[1)] The $k$-graph $G^*_k$ is asymmetric. 
		\item[2)] Every non-trivial $k$-subgraph of $G^*_k$ has an involution.
	\end{itemize}
\end{lemma} 	

\begin{proof}
	First, we prove that $G^*_k$ is asymmetric. Assume that $\phi$ is a non-trivial automorphism of $G^*_k$. $G_k$ is a $k$-subgraph of $G^*_k$, but in $G^*_k$, the degree of vertex $v_1$ is $2$ while $v_{2k-1}$ has degree $1$. Combining Lemma \ref{lem3}, we know that $\phi(v_i)= v_i$ for every $i\in [k]$. Then $\phi(x)=x$. Thus $G^*_k$ is asymmetric.
	
	To prove the second property of $G^*_k$,
	we assume $G$ is a non-trivial $k$-subgraph of $G^*_k$. Then either $G\subseteq G_k$ or $G$ is obtained by adding the vertex $x$ and the edge $\{x,v_1,\cdots,v_{k-2},v_{k+2}\}$ into a non-trivial $k$-subgraph of $G_k$. In the former case, $G$ has an involution by Lemma \ref{lem3}. 
	In the later case, let $i$ be the minimal index such that $G$ contains the edges $M_i, M_{i+1}, \dots, M_k$. Since $G$ is a nontrival $k$-subgraph of $G^*_k$, we have $i>1$ and $M_{i-1}$ is not in $G$. According to the proof of Lemma \ref{lem3}, we conclude $i\in\{2,3\}$. If $M_1$ is an edge of $G$, then $G$ has an involution $\phi$ such that $\phi(v_1)=v_2$ and $\phi(v_2)=v_1$ in the case $i=2$ (resp. $\phi(v_2)=v_3$ and $\phi(v_3)=v_2$ in the case $i=3$). Otherwise $M_1$ is not an edge of $G$, there is an involution $\phi$ such that $\phi(v_1)=x$ and $\phi(x)=v_1$. \qed
\end{proof}

For a hypergraph $G$, let $\widetilde{G}$ be a hypergraph obtained from $G$ by disjiontly adding two vertices to each edge of $G$.

\begin{observation}
	\label{obs1}
	For every hypergraph $G$, every automorphism of $\widetilde{G}$ which maps $G$ to $G$ is also an automorphism of $G$ and every automorphism of $G$ extends to $\widetilde{G}$.
\end{observation}

\begin{lemma}
	\label{lem5}
	Suppose $\phi$ is an automorphism of $\widetilde{G}_k=(\widetilde{X}_k,\widetilde{\mathscr{M}}_k)$ which leaves the set $\{v_1,v_{2k-2},v_{2k-1}\}$ invariant. Then $\phi$ induced on $G_k$ is identity.
\end{lemma}

\begin{proof}
	Observe that the degree of each vertex in $X_k\setminus\{v_1,v_{2k-2},v_{2k-1}\}$ is at least $2$ while every vertex in $\widetilde{X}_k\setminus X_k$ has degree $1$.
	As $\phi$ is an automorphism of $\widetilde{G}_k$ which leaves the set $\{v_1,v_{2k-2},v_{2k-1}\}$ invariant, $\phi$ maps $G_k$ to $G_k$.
	By Lemma \ref{lem3} and Observation \ref{obs1}, $\phi$ induced on $G_k$ is identity. \qed
\end{proof}

\begin{lemma}
	\label{lem6}
	Suppose $\phi$ is an automorphism of $\widetilde{G}^*_k=(\widetilde{X}^*_k,\widetilde{\mathscr{M}}^*_k)$ which leaves the vertices $x$ and $v_{2k-1}$ invariant. Then $\phi$ induced on $G^*_k$ is identity.
\end{lemma}

\begin{proof}
	The proof of this lemma is very similar to the above proof of Lemma \ref{lem5}.
	
	Observe that the degree of each vertex in $X^*_k\setminus\{x,v_{2k-1}\}$ is at least $2$ while every vertex in $\widetilde{X}^*_k\setminus X_k$ has degree $1$.
	As $\phi$ is an automorphism of $\widetilde{G}^*_k$ which leaves the set $\{x,v_{2k-1}\}$ invariant, $\phi$ maps $G^*_k$ to $G^*_k$.
	By Lemma \ref{lem4} and Observation \ref{obs1}, $\phi$ induced on $G^*_k$ is identity. \qed
\end{proof}

For $k\ge 6$ and any non-negative integer $s$, we shall construct a $k$-graph $G_{k,s}=(X,\mathscr{M})$.
Let $n=(k-1)(k-2)^s$.
First, we construct a hypergraph $H=(X,\hat{\mathscr{M}})$, depicted as Figure A3, which is consist of $s+2$ layers as follow: 
\begin{itemize}
	\item On layer 1, disjoint union of $n$ copies of $G_k$. 
	\item On layer 2, disjoint union of $\frac{n}{k-2}$ copies of $G_{k-2}$.
	\item On layer 3, disjoint union of $\frac{n}{(k-2)^2}$ copies of $G_{k-2}$.
	\item ...
	\item On layer $(s+1)$, disjoint union of $\frac{n}{(k-2)^s}=k-1$ copies of $G_{k-2}$.
	\item On layer $(s+2)$, one copy of $G^*_{k-2}$.
\end{itemize}

Intuitively, $G_{k,s}$ is obtained from $H$ by associating each $(k-2)$-edge in each copy of $G_{k-2}$ on Layer $(i+1)$ ($G^*_{k-2}$ on Layer $(s+2)$) to a copy of $G_{k-2}$ on Layer $i$ ($G_k$ on Layer $1$), $i\in [s+1]$ and changing each $(k-2)$-edge into a $k$-edge by adding the last two vertices of the corresponding copy of $G_{k-2}$ or $G_k$ to it.


\begin{figure}[!htp]
	\begin{center}
		\includegraphics[scale=0.2]{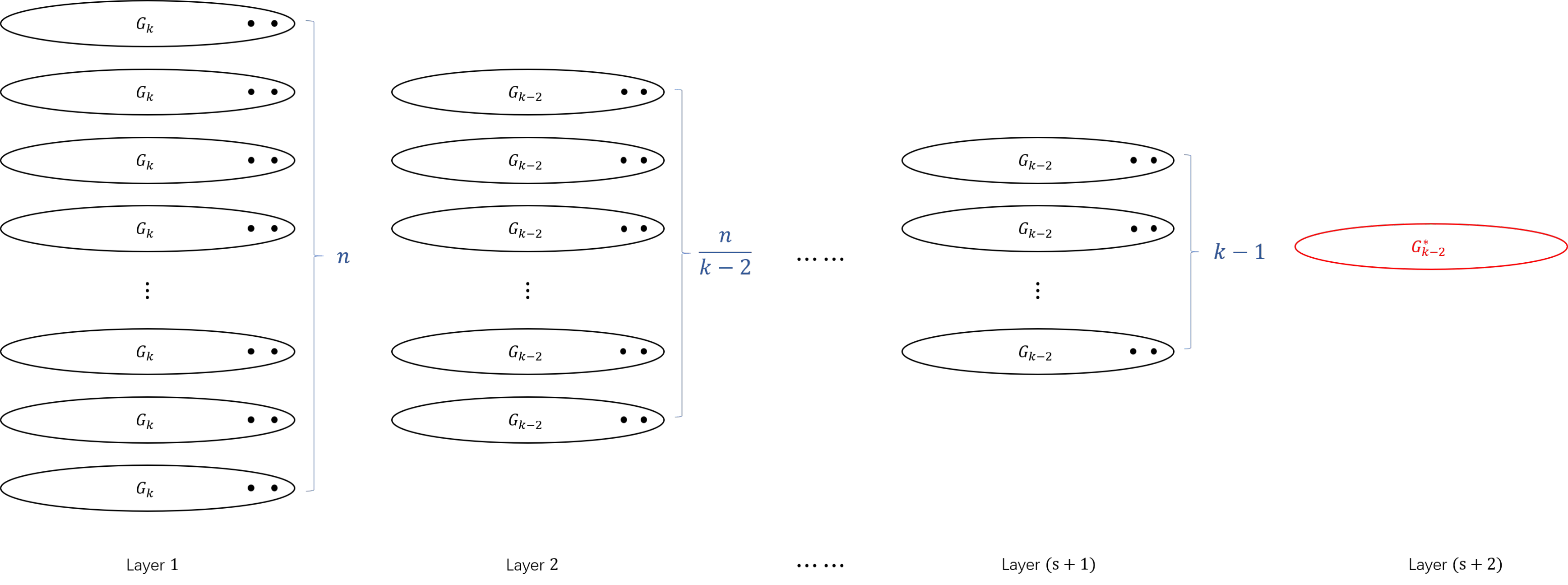}
	\end{center}
	\caption{The hypergraph $H$}
\end{figure}


Formally, this construction of $k$-graph $G_{k,s}=(X,\mathscr{M})$ can be done as follows:

The construction will proceed in two steps.

We shall put $n=(k-1)(k-2)^s$ and  we first consider $n$ copies of $G_k$, ${\frac{n}{(k-2)}}+{\frac{n}{(k-2)^2}}+\cdots+(k-1)$ copies of $G_{k-2}$ and one copy of $G^*_{k-2}$ arranged into $s+2$ layers (see schematic Figure A3). 

We have hypergraph $G^*_{k-2}$ on Layer $(s+2)$. Graphs on Layer $(s+1)$ are $k-1$ copies of $G_{k-2}$, which will be listed as $G(1)$, $G(2)$, $\cdots$, $G(k-1)$.
Graphs on Layer $s+1\le l\le 1$ will be $\frac{n}{(k-2)^{l-1}}$ copies of $G_{k-2}$ (or $G_k$ when $l=1$) and they will be listed as $G(i_l,i_{l+1},\cdots,i_{s+1})$, $1\le i_j\le k-2$, $j=l,l+1,\cdots,s$, $1\le i_{s+1}\le k-1$.

Then the vertices of $G_{k,s}$ are obtained from the vertices of the disjiont union of all hypergraphs $G(i_l,i_{l+1},\cdots,i_{s+1})$, $1\le l\le s+1$ and $G^*_{k-2}$.

Secondly we modify the $(k-2)$-edges to $k$-edges as follows:

Let $\mathscr{M}_{k-2}=\{M_1,M_2,\cdots,M_{k-2}\}$ and $\mathscr{M}^*_{k-2}=\{M^*_1,M^*_2,\cdots,M^*_{k-2},M^*_{k-1}\}$.
The edge $M^*_{i_l}$ of $G^*_{k-2}$ is enlarged by two last vertices of hypergraph $G(i_{s+1})$, $1\le i_{s+1}\le k-1$. Similarly, the edge $M$ correspongding to $M_{i_l}$ in $G(i_{l+1},\cdots,i_{s+1})$ on Layer $(l+1)$ is enlarged by the last two vertices of $G(i_l,i_{l+1},\cdots,i_{s+1})$ on Layer $l$, $1\le i_l\le k-2$. This finishes the construction of the $k$-graph $G_{k,s}$.

\medskip

Since $s$ can be any non-negative integer, it is sufficient to prove that for each of $k$-graphs $G_{k,s}$, the propertiess in Theorem \ref{main4} hold.

For the first property, to the contrary, we assume that $G_{k,s}$ has a non-trivial automorphism $\phi$.

\begin{claim}[A1]
	$\phi(v)=v$ for any vertex $v$ on Layer $(s+2)$,
\end{claim}

\begin{proof}
	Observing the degree sequence of each edge of $G_{k,s}$, the degree $1$  belong to four different types of degree sequences: the first edge in each copy $G(i_1,i_2,\cdots,i_{s+1})$ of $G_k$: $(1,2,\cdots,k)$; the corresponding edge $M$ of edge $M_{i_l}$ in each copy $G(i_l,i_{l+1},\cdots,i_{s+1})$ of $G_{k-2}$ on Layer $l$ ($2\le l\le s+2$): $(1,2,2,3,3\cdots,k-2)$; and the corresponding edges $M'$ and $M''$ of two different edges $M^*_1$ and $M^*_{k-1}$ of $G^*_{k-2}$ on Layer $(s+2)$: $(1,2,2,3,3\cdots,k-3,k-3)$ and $(1,2,2,3,3,\cdots,k-4,k-3,k-3,k-2)$.
	Thus the two vertices with degree $1$ on Layer $(s+2)$ are different from the others, which implies that $\phi$ maps the $\widetilde{G}^*_{k-2}$ to itself. But by Observation \ref{obs1} and  by Lemma \ref{lem6}, we obtain that $\phi$ restricted to Layer $(s+2)$ is identity. \qed
\end{proof}

\begin{claim}[A2]
	If the only automorphism $\phi$ of $G_{k,s}$ induced on Layer $(l+1)$ is identity, then $\phi$ induced on Layer $l$ is also identity, $l\ge 1$.
\end{claim}

\begin{proof}
	Since the automorphism $\phi$ of $G_{k,s}$ induced on Layer $(l+1)$ is identity, the corresponding edges $M$ of each edge $M_{i_l}$ in $G(i_{l+1},i_{l+2},\cdots,i_{s+1})$ on Layer $(l+1)$ are pairwise different. It implies that the copies $G(i_l,i_{l+1},\cdots,i_{s+1})$ on Layer $l$ are pairwise different and $\phi$ maps each $\widetilde{G}(i_l,i_{l+1},\cdots,i_{s+1})$ to itself (if $l=1$, $\phi$ maps each $G(i_1,i_2,\cdots,i_{s+1})$ to itself) and leaves the head vertex and the tail two vertices of $G(i_l,i_{l+1},\cdots,i_{s+1})$ invariant. 
	By Observation \ref{obs1},and by Lemma \ref{lem5}, the automorphism $\phi$ on $G$ induced on each copy $G(i_l,i_{l+1},\cdots,i_{s+1})$ of $G_{k-2}$ ($G_k$) on Layer $l$ is identity. \qed
\end{proof}

Claim (A1) states that the automorphism $\phi$ of $G_{k,s}$ induced on Layer $(s+2)$ is identity. Then by Claim (A2), $\phi$ of $G_{k,s}$ induced on Layer $(s+1)$ is identity. Continuing this way, we obtain that $\phi$ induced on Layer $i$ is identity, $i\in [s+2]$. Thus $G_{k,s}$ is asymmetric.

For the second property, we have the followig claim.

\begin{claim}[A3]
	For every $k\ge 6$ and $s\ge 1$, any proper $k$-subgraph of $G_{k,s}$ with at least $2$ vertices has an involution.
\end{claim}

\begin{proof}
	For contradiction, assume that $G_{k,s}$ contains a non-trivial $k$-subgraph $H$ such that $H$ is asymmetric. Without loss of generality, let us assume that $H$ is connected.
	
	Now if $H$ fails to contain all vertices and edges on Layer $1$, then there exists a copy $G=G(i_1,i_2,\cdots,i_{s+1})$ of $G_k$ such that $G'=H\cap G$ is a proper $k$-subgraph of $G$.
	
	If $G'$ has at least $2$ vertices, then by Lemma \ref{lem3}, $G'$ has an involution $\phi$ which leaves the tail two vertices $x$, $y$ of $G$ invariant in the case that $\{x,y\}\subsetneqq V(G')$. When $G'=\{x,y\}$, it is obvious that $G'$ has an involution. In the case that one of the vertices $x$, $y$ is missing in $G'$, $G'$ is not connected to the rest of $H$. If $G'$ is empty then let $l$ be the minimal $l$ such that the hypergraph $G''$ induced on $G(i_l,i_{l+1},\cdots,i_{s+1})$ by $H$ is non-empty.
	Then for the minimality of $l$ and the connectivity of $H$, either $G''$ is exactly the tail two vertices of $G(i_l,i_{l+1},\cdots,i_{s+1})$ thus $G''$ has an involution $\phi$, or $G''$ is empty. If $G''$ is empty, then assume that $H$ contains some vertices of $G^m=G(j_m,j_{m+1},\cdots,j_{s+1})$, $j_t\in [k-2], m\le t\le s+1$. Let $G'''=H\cap G^m$. Then by Observation \ref{obs1}, there is a non-identical automorphism (an involution) $\phi$ of $\widetilde{G}'''$ which maps $G'''$ to $G'''$. In all above cases, $\phi$ extends to whole $H$. 
	
	We conclude that the whole Layer $1$ has to be present in $H$. It then follows that all edges connecting Layer $1$ to Layer $2$ have to be present in $H$ for otherwise $H$ would not be connected. But this then implies that all edges on Layer $1$ and $2$ are present in $H$. Continuing this way, we obtain that $H$ contains all edges of $G_{k,s}$, a contradiction. \qed
\end{proof}	

\medskip
Therefore our proof is complete.
And it is easy to observe that the $k$-graphs $G_{k,s}$ have vertex degree bounded by $k$.

\end{document}